\documentclass[11pt]{article}

\usepackage[margin = 0.8in]{geometry}
\usepackage{amssymb}
\usepackage{amsmath}
\usepackage[utf8]{inputenc}
\usepackage[T1]{fontenc}   
\usepackage{hyperref}      
\usepackage{url}           
\usepackage{booktabs}      
\usepackage{amsfonts}      
\usepackage{nicefrac}      
\usepackage{microtype}     
\usepackage{xcolor}        
\usepackage{amsmath, amsthm}
\usepackage{graphicx}
\usepackage{multirow}
\usepackage{wrapfig}
\usepackage{natbib} 
\usepackage{appendix}
\usepackage{comment}

\newtheorem{theorem}{Theorem}[section]

\theoremstyle{definition}

\newtheorem{remark}[theorem]{Remark}

\newcommand{\R}{\mathbb{R}}

\newcommand{\xx}{\mathbf{x}}
\newcommand{\aaa}{\boldsymbol{a}}
\newcommand{\dx}{\,\mathrm{d}\mathbf{x}}

\newcommand{\NN}{\mathcal{N}}
\newcommand{\LL}{\mathcal{L}}
\newcommand{\om}{\boldsymbol{\omega}}

\DeclareMathOperator*{\argmin}{arg\,min}

\begin{document}

\title{Data-Free Asymptotics-Informed Operator Networks for Singularly Perturbed PDEs}

\author{
Jinsil Lee\thanks{Department of Mathematical Sciences, Seoul National University, Seoul, Republic of Korea. Email: \tt{jl74942@snu.ac.kr}},
~Youngjoon Hong\thanks{Department of Mathematical Sciences \& Research Institute of Mathematics, Seoul National University, Seoul, Republic of Korea. Email: \tt{hongyj@snu.ac.kr}},
~Seungchan Ko\thanks{Department of Mathematics, Inha University, Incheon, Republic of Korea. Email: \tt{scko@inha.ac.kr}},
~and~Jae Yong Lee\thanks{Department of AI, Chung-Ang University, Seoul, Republic of Korea. Email: \tt{jaeyong@cau.ac.kr}}
}

\date{}
\maketitle

\begin{abstract}
Recent advances in machine learning (ML) have opened new possibilities for solving partial differential equations (PDEs), yet robust performance in challenging regimes remains limited. In particular, singularly perturbed differential equations exhibit sharp boundary or interior layers with rapid transitions, where standard ML surrogates often fail without extensive resolution. Generating training data for such problems is also costly, as accurate reference solutions typically require massive adaptive mesh refinement.
In this work, we propose eFEONet, an enriched Finite Element Operator Network tailored to singularly perturbed problems. Guided by classical singular perturbation theory, eFEONet augments the operator-learning framework with specialized enrichment basis functions that encode the asymptotic structure of layer solutions. This design enables accurate approximation of sharp transitions without relying on large datasets, and can operate with minimal supervision, or even in a data-free manner under appropriate settings. We further provide a rigorous convergence analysis of the proposed method and demonstrate its effectiveness through extensive experiments on representative problems featuring both boundary and interior layers.
\end{abstract}

\section{Introduction}\label{sec1}
\begin{figure}[t]
   \centering
    \includegraphics[width=.8\textwidth]{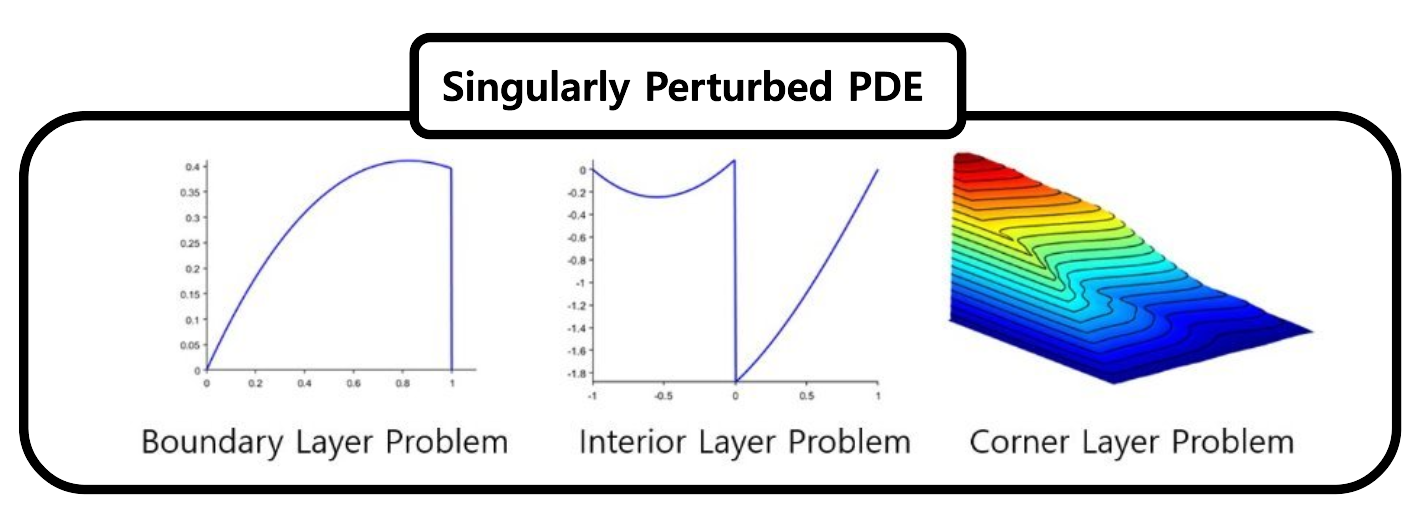}
   \caption{Representative solution profiles for singularly perturbed PDEs, illustrating the inherent stiffness of boundary and interior layers across various domains. The sharp gradients and rapid transitions depicted here highlight the intrinsic stiffness and associated computational challenges.}
   \label{fig:singularyPDE}
   \vspace{-3mm}
\end{figure}

The use of machine learning (ML) to solve partial differential equations (PDEs) has made significant advancements in recent years, offering innovative approaches to tackle longstanding challenges in scientific computing \citep{lagaris1998artificial, PINN001, yu2018deep,ainsworth2021galerkin}. Among these methods, operator networks have emerged as a practical and efficient tool due to their ability to infer solutions quickly after training \citep{lu2021learning,li2021fourier}. Unlike classical numerical methods, which solve PDEs iteratively for each new set of conditions, operator networks enable rapid predictions by learning the solution operator itself. This capability has made operator networks a promising new paradigm in the study of parametric PDEs. However, there remain notable limitations in applying operator networks to real-world problems. One key challenge is that training an operator network typically requires a pre-generated dataset of solutions, which are often constructed using conventional numerical methods. This process can be computationally expensive, particularly for complex PDEs. Singularly perturbed differential equations, in particular, present unique difficulties. These equations often exhibit rapid transitions within thin regions known as boundary or interior layers (see Figure~\ref{fig:singularyPDE}), where constructing accurate datasets is both costly and technically challenging. Moreover, the sharp transitions inherent in thin layers can degrade the performance of operator networks, which tend to rely on smooth priors \citep{lu2022comprehensive}.

\begin{figure}[t]
   \centering
       \includegraphics[width=0.42\textwidth]{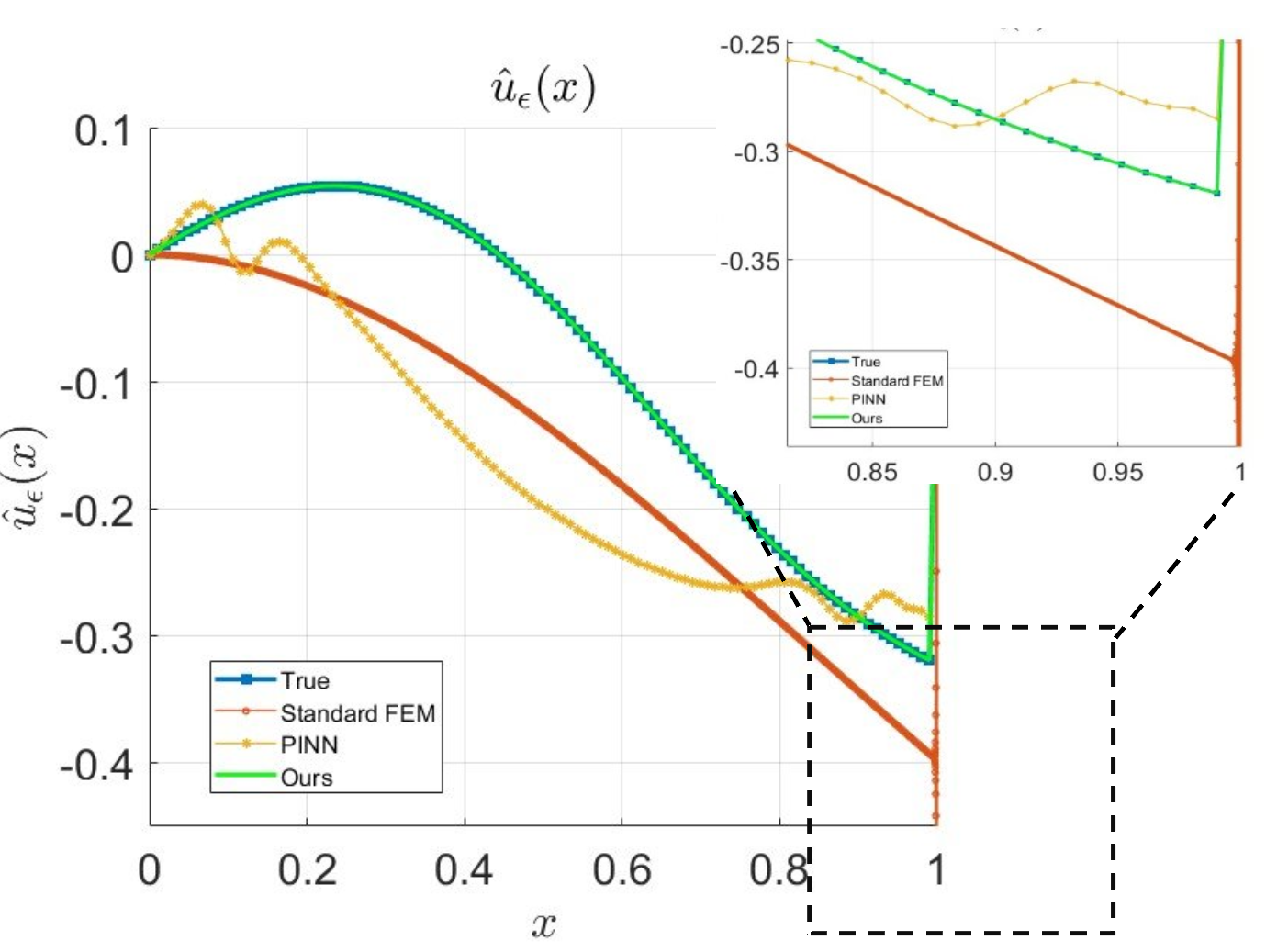}
    \includegraphics[width=0.42\textwidth]{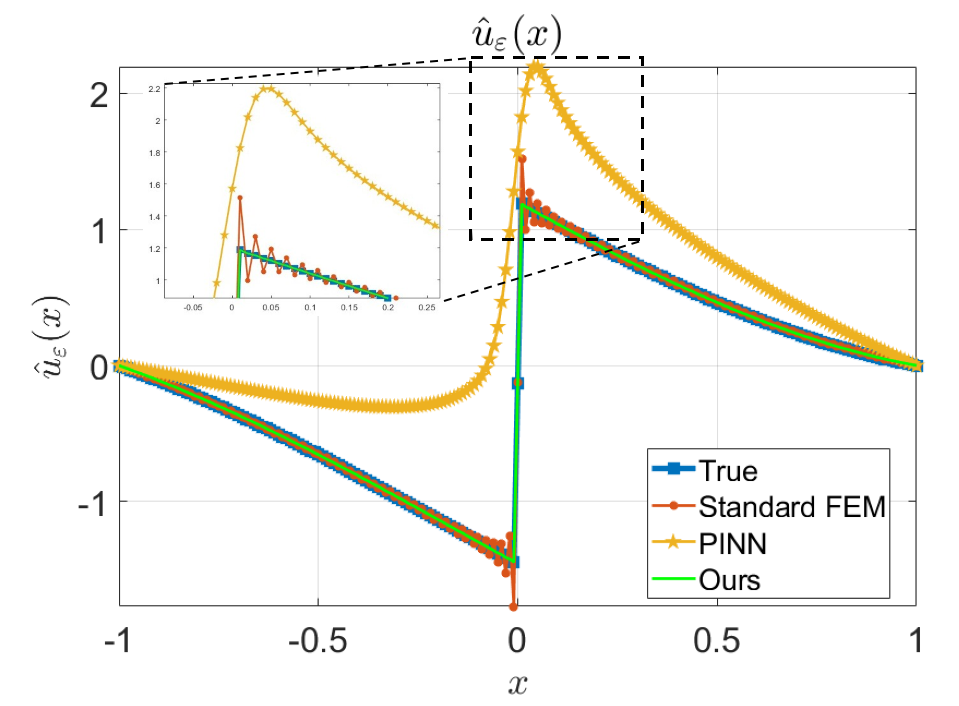}
   \caption{Examples of the boundary (left) and interior (right) layer phenomena and comparisons of the reference solution (True) and the predicted solutions using Standard FEM, PINN, and eFEONet (Ours). We set $\varepsilon = 10^{-5}$, the mesh size for the Standard FEM is $1/12000$ for the left case, while it is kept identical to those of the other methods for the right case.}
   \label{fig:ex_layers}
   \vspace{-3mm}
\end{figure}

Boundary and interior layer phenomena are of paramount importance in many scientific and engineering disciplines, including fluid dynamics, biology, and chemical reactions \citep{schlichting2016boundary, batchelor2000introduction}. These problems are characterized by sharp changes in solution profiles within thin layers, making them notoriously difficult to handle even with advanced numerical methods. The challenge arises from the small diffusive parameter $\varepsilon>0$ in these equations, which leads to steep gradients over small spatial regions. See Figure \ref{fig:ex_layers} where the examples of boundary and interior layer phenomena are presented. Developing methods to accurately and efficiently solve such problems remains a challenging task in scientific computing \citep{FEM1,FEM2}. 
Machine learning-based approaches face additional challenges because they are inherently better at learning smooth functions but struggle to accurately capture sharp transitions. Neural networks, for instance, are often designed to approximate solutions that vary gradually, making it difficult to capture the steep gradients and sharp transitions characteristic of boundary or interior layers \citep{PINN007}.
These boundary and interior layer problems require computationally expensive massive mesh refinement to obtain accurate solutions. Moreover, as \( \varepsilon>0\) decreases, the mesh size must become finer, following an approximate scaling of adaptive mesh size $\simeq \varepsilon$. 
This results in a significant drawback, as data generation can become prohibitively expensive in many cases.
This limitation highlights the need for new architectures and methodologies that can handle these complexities without compromising accuracy or efficiency.

Operator learning trains models to approximate PDE solution operators using datasets of input-output pairs from numerical solvers \citep{bhatnagar2019prediction, guo2016convolutional, khoo2017solving, zhu2018bayesian}, enabling efficient and real-time predictions for varying inputs \citep{NEURIPS2020_4b21cf96}. Notable architectures include the Fourier Neural Operator (FNO) \citep{kovachki2021neural} and DeepONet \citep{lu2021learning}. Recent advances also explore message-passing frameworks to accommodate complex problem structures \citep{brandstetter2022message, lienen2022learning, pfaff2021learning, boussif2022magnet}. In addition, transformer-based architectures have been introduced \citep{cao2021choose,wang2025cvit,pmlr-v235-hao24d}, along with emerging foundation models tailored for PDEs\citep{herde2024poseidon,ye2024pdeformer}. Despite these developments, operator learning still faces challenges in generalization, data efficiency, and resolving sharp solution features. Among various operator-learning models, unsupervised physics-based operator networks incorporate governing equations directly into neural operator architectures, minimizing or completely removing the need for labeled training data. Variational frameworks such as Finite Element Operator Network (FEONet) \citep{lee2024finite} and Spectral Coefficient Learning via Operator Network (SCLON) \citep{SCLON} use PDE residuals in weak form to achieve accurate predictions without explicit simulation data. Similarly, physics-informed approaches like Physics-Informed Neural Networks (PINNs) \citep{PINN001, han2018solving}, PINO \citep{PINO}, and PIDeepONet \citep{PIDON} can also be formulated to rely entirely on PDE constraints and boundary conditions. Despite recent progress, accurately capturing multiscale phenomena and sharp gradients without labeled data remains challenging, highlighting the need for more robust unsupervised approaches.

In recent years, deep learning has emerged as a promising approach for solving singularly perturbed PDEs, with physics-informed methods also contributing to this line of research \citep{SP01, SP02}. However, many of these approaches suffer from limited scalability and tend to be effective only in restricted settings. For instance, the study on stiff chemical kinetics in \citep{GOSWAMI2024116674} employs deep neural operators tailored to reaction-diffusion stiffness, which narrows its applicability. In contrast, our method targets a broader class of singularly perturbed PDEs, including problems with boundary and interior layers, and is particularly suited to data-scarce regimes. More recently, \citep{chen2025learn} proposed a homotopy-based learning strategy for singularly perturbed problems, but it focuses on specific PDE instances rather than an operator-learning framework. Component Fourier Neural Operator (ComFNO) \citep{SP03}, a modified variant of FNO, also attempts to incorporate asymptotic behavior into the model. Nonetheless, key challenges remain, including the need for large training datasets, difficulties in accurately capturing sharp transitions, and the absence of rigorous theoretical foundations to ensure reliability across a wide range of problems.

In this paper, we propose an enriched Finite Element Operator Network (eFEONet), specifically designed to address these challenges. eFEONet builds upon the FEONet framework \citep{lee2024finite}, a highly data-efficient operator learning method that requires minimal training data or no dataset at all.
Unlike traditional operator networks, eFEONet leverages the structure of finite element methods (FEMs), where the solution is expressed as a linear combination of nodal coefficients and basis functions. This design not only eliminates the need for large datasets but also ensures the exact satisfaction of boundary conditions.  
By incorporating insights from well-developed singular perturbation analysis in PDE theory, we design special basis functions within the finite element framework that capture the asymptotic behavior of solutions in boundary or interior layers \citep{BLA_3}. This approach enables accurate modeling of sharp transitions while maintaining computational efficiency. Furthermore, as demonstrated in Section~\ref{sec:experiments}, eFEONet significantly outperforms data-driven neural operator baselines in the singularly perturbed regime. In particular, while ComFNO is trained using 900 input-output samples, it exhibits errors that are two orders of magnitude larger than those of eFEONet, which is trained without any paired solution data.
We validate our approach through theoretical convergence analysis and empirical results on various singularly perturbed problems, including both boundary and interior layers \citep{H2020, BLA_3}. The results demonstrate that eFEONet achieves high accuracy and efficiency, even for problems with strong boundary or interior layer phenomena such as convection-dominated PDEs \citep{Stynes_2005}.

The main contributions of the paper are summarized as follows: First, we propose eFEONet, which integrates singular perturbation analysis into the FEONet framework. This incorporation enables superior accuracy in solving singularly perturbed PDEs, effectively capturing sharp transitions in both boundary and interior layers.  Secondly, Singular perturbation problems typically require increasingly finer meshes as the parameter \( \varepsilon>0 \) decreases, making dataset generation computationally very expensive. Our approach overcomes this limitation by being highly data-efficient, requiring minimal training data, or even operating without any training dataset. Finally, we demonstrate the effectiveness of eFEONet through comprehensive experiments on challenging convection-diffusion PDEs, including problems with boundary and interior layers in both 1D and 2D. The results show that eFEONet achieves error reductions of two orders of magnitude compared to existing approaches, even when no training data is used.  

\section{Proposed Method: Enriched Finite Element Operator Network}

In this section, we shall describe our proposed method, eFEONet, designed for solving singularly perturbed parametric PDEs. We start by giving a brief overview of the enriched FEMs, which form the core of our approach. Then, we will explain eFEONet, the main method we propose in this paper.

For the description, we will focus on the following PDE:
\begin{equation}\label{main_eq}
\begin{aligned}
-\varepsilon\,{\rm{div}}\,(\aaa(\xx)\nabla u_{\varepsilon})+\boldsymbol{b}(\xx)\cdot\nabla u_{\varepsilon}&=f\quad{\rm{in}}\,\,D,\\
u_{\varepsilon}&=0\quad{\rm on}\,\,\partial D.
\end{aligned}
\end{equation}
Here we assume that the singular perturbation parameter $\varepsilon>0$ is very small so that the boundary layer phenomenon occurs. Furthermore, to highlight that the shape of a solution depends on $\varepsilon>0$, we will denote the solution as $u_{\varepsilon}$. As we will explain in more detail later, we propose an operator-learning approach for the singular perturbation problem that enables real-time solution predictions whenever the input data of the PDE varies. As a prototype model, we set the external force  $f$ as an input of neural networks, and train the model so that the neural networks can learn the operator $\mathcal{G}:f\mapsto u_{\varepsilon}$. Note, however, that our method can be easily extended to various forms of input functions, including boundary conditions, variable coefficients, or initial conditions (see, e.g., \citep{lee2024finite}).

\subsection{Finite Element Method (FEM)}\label{desc_FEM}

The FEM is a general technique for the numerical solution of PDEs. It is based on the variational formulation of the PDE \eqref{main_eq}, which seeks to find a function $u_{\varepsilon}\in V$ satisfying
\begin{equation}\label{var_for}
B[u_{\varepsilon},v]:=\varepsilon\int_{D}\aaa(\xx)\nabla u_{\varepsilon}\cdot\nabla v\dx+\int_{D}\boldsymbol{b}(\xx)\cdot\nabla u_{\varepsilon}v\dx=\int_{D}fv\dx=:\ell(v)\quad{\text{for all}}\,\,v\in V,
\end{equation}
where $V$ is typically an infinite-dimensional function space for the solution and test functions. The first step in FEM theory is to discretize the domain \( D \subset \mathbb{R}^d \), known as a \textit{triangulation}. For \( d = 1 \) and \( D = [a, b] \), this involves points $
a = x_0 < x_1 < \cdots < x_K = b,$
with each interval \( [x_{i-1}, x_i] \) forming a 1-simplex. For \( d = 2 \), the triangulation consists of closed triangles \( T_i \) (2-simplexes), \( i = 1, \ldots, K \), whose interiors are disjoint. If \( i \neq j \) and \( T_i \cap T_j \neq \emptyset \), then the intersection is either a shared vertex or edge. For \( d \geq 3 \), elements are \( d \)-simplexes. Let \( h_T \) denote the longest edge of a triangle \( T \), and define the global mesh size as $h = \max_T h_T$. Let \( S_h \) be the space of continuous functions \( v_h \) on \( D \) such that the restriction of \( v_h \) to each element is a polynomial. The finite-dimensional ansatz space is then defined as $V_h = S_h \cap V$. Let \( \{ \mathbf{x}_i \} \) denote the triangulation vertices, and \( \{ \phi_j \} \) the \textit{nodal basis} for \( V_h \) with $\phi_j(\mathbf{x}_i) = \delta_{ij}$.
Using piecewise linear basis functions defines the P1-element method; using piecewise quadratic polynomials gives the P2-element method. The dimension of \( V_h \) depends on the triangulation and hence on the mesh parameter \( h \).

The FEM aims to approximate the infinite-dimensional space $V$ by a finite-dimensional subspace $V_h$ defined by $V_h={\rm span}\{\phi_1,\phi_2,\cdots,\phi_{N(h)}\}
$.
This makes the problem numerically solvable. Motivated from \eqref{var_for}, we seek to compute the approximate solution $u_{\varepsilon,h}\in V_h$ using the so-called {\textit{Galerkin approximation}}, which is given by the equation
\begin{equation}\label{gal_approx}
B[u_{\varepsilon,h},v_h]=\ell(v_h)\quad{\text{for all}}\,\,v_h\in V_h.
\end{equation}
Writing the finite element solution as
$u_{\varepsilon,h}(\xx)=\sum^{N(h)}_{k=1}\alpha_k\phi_k(\xx)$ with $\alpha_i\in\R$, the Galerkin approximation \eqref{gal_approx} transforms into the following linear algebraic system:
\begin{equation}\label{lin_alg_eq}
A\alpha=F\quad{\text{with}}\,\,A_{ik}:=B[\phi_k,\phi_i]\,\,{\rm{and}}\,\,F_i:=\ell(\phi_i).
\end{equation}
Here, the matrix $A$ is invertible, assuming the underlying PDE has an appropriate structure. The coefficients $\{\alpha_k\}^{N(h)}_{k=1}$ can be determined by solving \eqref{lin_alg_eq}, thus yielding the approximate solution $u_{\varepsilon,h}$.

\subsection{Enriched FEONet with a Corrector Basis}
\begin{figure}[t]
   \centering
    \includegraphics[width=1.0\textwidth]{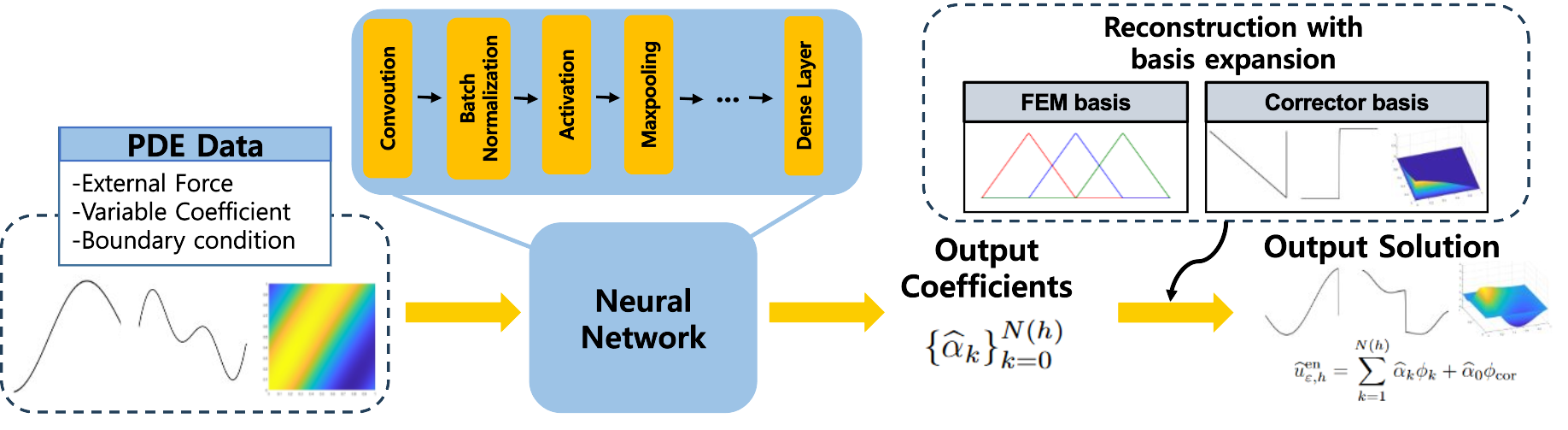}
   \caption{Schematic illustration of eFEONet.}
   \label{fig:scheme}
\end{figure}
In this section, we shall introduce our main method, the enriched FEONet (eFEONet). One key novelty of the eFEONet is to utilize extra basis functions derived from theoretical arguments in \citep{BLA_3}. For a clear illustration of the proposed method, we shall explain it through a simple example of the following form:
\begin{equation}\label{eq_singular}
  \begin{aligned}
    -\varepsilon u''_{\varepsilon}-u'_\varepsilon&=f(x), \quad x \in (-1,1),\\
    u_\varepsilon(-1)&=u_\varepsilon(1)=0, 
  \end{aligned} 
\end{equation}
where $0<\varepsilon \ll 1$.
As shown in Figure \ref{fig:ex_layers}, when $\varepsilon>0$ is small, it is difficult to expect other well-known methods, including classical FEM and PINN, to achieve good performance due to the sharp transitions near the boundary.
To accurately capture the boundary layer, we incorporate an additional basis function, commonly referred to as the {\textit{corrector function}} in mathematical analysis, for example in this case, defined as: $\phi_{\rm cor}(x):=e^{-(1+x)/\varepsilon}-(1-(1-e^{-2/\varepsilon})(x+1)/2)$. Such a basis function reflects the boundary layer properties of the given equation and is derived from theoretical arguments. The derivation of various corrector basis functions can be found in a range of references, including \citep{BLA_3}. The corrector basis is added to the standard nodal basis functions of FEM to construct an enriched Galerkin space. In other words, for enriched FEONet for the singularly-perturbed problems, we now replace the original ansatz space $V_h$ by the enriched Galerkin space $\overline{V}_h:=\{\phi_{\rm cor}, \phi_1,\phi_2,\cdots,\phi_{N(h)}\}$,
where the corrector basis $\phi_{\rm cor}$ has been added to $V_h$. It is noteworthy that no high additional computational cost occurs, as the enriched basis is only restricted to boundary elements. In general, neural networks assume a smooth prior, which makes them less effective in handling boundary layers. This can lead to unstable training due to the direct calculation of the PDE residual. In contrast, the eFEONet leverages theory-guided basis functions, allowing its predicted solution to capture the sharp transitions near the boundary precisely. Encapsulating the above discussion, the enriched FEM for the boundary layer problem can be written as follows: we seek $u^{\rm en}_{\varepsilon,h}\in \overline{V}_h$ satisfying
\begin{align}\label{en_gal_approx}
B[u^{\rm en}_{\varepsilon,h},v_h]:&=\varepsilon\int_{D}\aaa(\xx)\nabla u^{\rm en}_{\varepsilon,h} \cdot\nabla v_h\dx+\int_{D}\boldsymbol{b}(\xx)\cdot\nabla u^{\rm en}_{\varepsilon,h} v_h\dx \notag\\
&=\int_{D}fv_h\dx=:\ell(v_h)\quad{\text{for all}}\,\,v_h\in \overline{V}_h.
\end{align}
In our eFEONet approach, the input to the neural network consists of data related to the given PDE problems, which is parameterized by $\om\in\Omega$, while the output consists of the coefficients of a basis expansion. To be more specific, we incorporate this into a deep learning framework to construct the eFEONet, whose solution prediction is written as
\begin{equation}\label{sol_recon_new}
    \widehat{u}^{\rm en}_{\varepsilon,h}(\xx;\om)=\sum_{k=1}^{N(h)}\widehat{\alpha}_k(\om)\phi_k(\xx)+\widehat{\alpha}_0(\om)\phi_{\rm cor},
\end{equation}
where the dimension of the output of the neural network has increased by one to handle the added corrector basis.
By writing $\phi_0:=\phi_{\rm cor}$, the loss function for the eFEONet is defined as
\begin{equation}\label{new_loss_function}
\LL^M(\widehat{u}^{\rm en}_{\varepsilon,h})=\frac{1}{M}\sum^M_{m=1}\sum^{N(h)}_{i=0}|B[\widehat{u}^{\rm en}_{\varepsilon,h}(\xx;\om_m),\phi_i(\xx)]-\ell(\phi_i(\xx);\om_m)|^2,
\end{equation}
for randomly drawn parameters $\om_1,\cdots,\om_M\in\Omega$. A schematic diagram of the eFEONet algorithm is shown in Figure \ref{fig:scheme}.

\begin{remark}
Our framework employs corrector functions tailored to specific problem classes, yet they are not confined to individual instances. For families of PDEs with analogous singular behavior, the same correctors can often be applied effectively. In convection-diffusion equations, for example, the boundary layer typically has a thickness proportional to $\varepsilon$ with an exponential profile, a structure preserved even with additional reaction terms.
\end{remark}

\begin{remark}
Some preliminary results show that one could attempt to learn the corrector bases using data. In contrast, our approach constructs them via classical numerical analysis, which not only requires no data but also achieves substantially better performance. This integration of analytical methods into an operator-learning framework constitutes the main novelty of our work, highlighting how analytic knowledge can maximize the efficiency of operator learning.
\end{remark}
\subsection{Convergence of eFEONet}
In this section, we discuss the convergence result for eFEONet, providing a theoretical foundation for the proposed approach. Although the convergence analysis of eFEONet does not differ substantially from that of the original FEONet \cite{lee2024finite}, we include a concise exposition in this section in order to introduce a coherent framework for the theoretical analysis of methods of this type. We let an external forcing term $f$ as the input of neural networks, that is parametrized by $\om$ in the probability space $(\Omega,\mathcal{F},\mathbb{P})$. We shall interpret $f(\xx;\om)$ as a bivariate function defined on $D\times\Omega$. Moreover, we will assume that
\begin{equation}\label{f_ass}
    f(\xx;\om)\in C(\Omega;L^1(D)):=\left\{f:\Omega\rightarrow L^1(D):\sup_{\om\in\Omega}\int_D|f(\xx;\om)|\dx<\infty\right\}.
\end{equation}
For each $\om\in\Omega$, the external force $f(\xx;\om)$ is specified, and the corresponding weak solution is denoted by $u_{\varepsilon}(\xx;\om)$. For given mesh size $h>0$, let $\overline{V}_h\subset H^1_0(D)$ be a finite-dimensional space spanned by the basis functions $\{\phi_k\}_{k=0}^{N(h)}$ including the corrector basis function $\phi_0=\phi_{\rm cor}$,
and $u^{\rm en}_{\varepsilon,h}\in \overline{V}_h$ be an enriched finite element approximation of $u_\varepsilon$ which satisfies the enriched Galerkin approximation \eqref{en_gal_approx}.
We write 
\begin{equation}\label{FEM_sol}
    u^{\rm en}_{\varepsilon,h}(\xx,\om)=\sum^{N(h)}_{k=0}\alpha^*_k(\om)\phi_k(\xx),
\end{equation}
where $\alpha^*$ is the finite element coefficients obtained from solving the corresponding linear algebraic system. Note that $\alpha^*$ can also be characterized in an alternative way:
\begin{equation}
    \alpha^*=\argmin_{\alpha\in C(\Omega,\R^{N(h)+1})}\LL(\alpha),
\end{equation}
where $\LL$ is the population risk
\begin{equation}\label{anal_pop_loss}
    \LL(\alpha)=\mathbb{E}_{\om\sim\mathbb{P}_{\Omega}}\bigg[\sum^{N(h)}_{i=0}|B[\widehat{u}(\om),\phi_i]-\ell(\phi_i;(\om))|^2\bigg].
\end{equation}

Next, we define the class of feed-forward neural networks as $\NN_n$, where the subscript $n$ denotes the network architecture. We assume that $\NN_{n_2}$ is more expressive than $\NN_{n_1}$ when $n_1 \leq n_2$. For instance, $n$ could represent the number of layers with bounded width, or the number of neurons when the number of layers is fixed. Neural networks are known to be an appropriate choice for nonlinear approximation, supported by the universal approximation theorem (see, for example, \citep{UA_1, UA_2, UA_3, terry_lion}). Now for a neural-network approximation of $\alpha^*$, we mean that $\widehat{\alpha}(n):\Omega\rightarrow\R^{N(h)+1}$, which solves the following minimization problem
\begin{equation}\label{for_2}
    \widehat{\alpha}(n)=\argmin_{\alpha\in\NN_n}\LL(\alpha),
\end{equation}
and we write the corresponding solution prediction by
\begin{equation}\label{for_3}
    \widehat{u}^{\rm en}_{\varepsilon,h,n}(\xx;\om)=\sum^{N(h)}_{k=0}\widehat{\alpha}(n)_k(\om)\phi_k(\xx).
\end{equation}
Note here that for the neural network $\alpha\in\NN_n$, the input is $\om\in\Omega$ that specifies the external forcing term $f(\xx;\om)$ and the output is the coefficient vector in $\R^{N(h)+1}$.

Finally, we define the solution of the following discrete minimization problem:
\begin{equation}\label{for_4}
    \widehat{\alpha}(n,M)=\argmin_{\alpha\in\NN_n}\LL^M(\alpha).
\end{equation}
Here $\LL^M$ is the empirical risk, which is the Monte Carlo integration of the population risk $\LL(\alpha)$:
\begin{equation}\label{anal_emp_loss}
    \LL^M(\alpha)=\frac{|\Omega|}{M}\sum_{m=1}^M\sum^{N(h)}_{i=0}|B[\widehat{u}(\om_m),\phi_i]-\ell(\phi_i;(\om_m))|^2,
\end{equation}
where $\{\om_n\}^M_{m=1}$ is an i.i.d. random variables following $\mathbb{P}_{\Omega}$.
 We then write the associated solution as
\begin{equation}\label{for_5}
    \widehat{u}^{\rm en}_{\varepsilon, h,n,M}(\xx;\om)=\sum^{N(h)}_{k=0}\widehat{\alpha}(n,M)_k(\om)\phi_k(\xx),
\end{equation}
which is the actual solution prediction by eFEONet. We assume that we can always find the exact minimizers for the problems \eqref{for_2} and \eqref{for_4}, and the optimization error is ignorable.

Let us denote the solution of \eqref{eq_singular} by $u_{\varepsilon}$ corresponding to a given parameter $0<\varepsilon\ll1$. 
Since our method is built upon the enriched FEM, the enriched finite element approximation $u^{\rm en}_{\varepsilon,h}$ in \eqref{en_gal_approx} serves as an intermediate step between the exact solution $u_{\varepsilon}$ and the approximate solution $\widehat{u}^{\rm en}_{\varepsilon,h,n,M}$ obtained from eFEONet \eqref{sol_recon_new}. Specifically, the total error $u_{\varepsilon}-\widehat{u}^{\rm en}_{\varepsilon,h,n,M}$ is decomposed into three parts:
\begin{equation}\label{error_split}
    u_{\varepsilon}-\widehat{u}^{\rm en}_{\varepsilon,h,n,M}=
    (u_\varepsilon-u^{\rm en}_{\varepsilon,h})+
    (u^{\rm en}_{\varepsilon,h}-\widehat{u}^{\rm en}_{\varepsilon,h,n})+
    (\widehat{u}^{\rm en}_{\varepsilon,h,n}-\widehat{u}^{\rm en}_{\varepsilon,h,n,M}).
\end{equation}
The first error arises from the finite element approximation, which we assume to be negligible when $h>0$ is sufficiently small.
The error analysis for the first term is well investigated in the previous literature on singular perturbation analysis. For example, in \citep{BLA_1}, the following error estimate was derived for the enriched FEM \eqref{en_gal_approx}:
\begin{equation}\label{egal_error}
    \|u_{\varepsilon}-u^{\rm en}_{\varepsilon,h}\|_{H^1}\leq C\left(h+\frac{h^2}{\varepsilon}\right),
\end{equation}
where $C>0$ is a constant independent of $h$ and $\varepsilon$. This result is especially highlighted as it provides a satisfactory convergence result even in the under-resolved case for $h>\varepsilon$. A mathematical analysis of this problem constitutes an independent topic traditionally addressed within classical numerical analysis; accordingly, we shall not examine it in detail in the present work. More general results can be found in various papers, e.g., from \citep{BLA_2, BLA_3}. For the analysis of eFEONet, based on the estimate \eqref{egal_error}, we can reduce this error to any desired level by selecting a suitable $h>0$. Therefore, we assume that $h$ has been chosen so that the finite element approximation error is small enough. The second error, known as the {\textit{approximation error}}, occurs when we use a class of neural networks to approximate the target (finite element) coefficients. The third error, often referred to as the {\textit{generalization error}}, measures how well our approximation performs on unseen data. Our focus will be on proving that, with fixed $h>0$ and $\varepsilon>0$, as the index $n\in\mathbb{N}$ for neural network architectures becomes larger and the number of input samples $M\in\mathbb{N}$ increases, our approximate solution $\widehat{u}^{\rm en}_{\varepsilon,h,n,M}$ converges to the finite element solution $u^{\rm en}_{\varepsilon,h}$ which is assumed to be the true solution here.

\begin{theorem}[Convergence of eFEONet]\label{main_thm_whole}
Assume that \eqref{f_ass} holds. Then for given $\varepsilon>0$ and $h>0$, with probability $1$, we have that
\begin{equation}\label{main_conv_whole}
    \lim_{n\rightarrow \infty}\lim_{M\rightarrow \infty}\|u^{\rm en}_{\varepsilon,h}-\widehat{u}^{\rm en}_{\varepsilon,h,n,M}\|_{L^2(\Omega;L^2(D))}=0.
\end{equation}
\end{theorem}
As mentioned earlier, this theorem can be shown directly from the proof presented in the original paper concerning FEONet \cite{lee2024finite}. The only difference lies in the singular perturbation analysis, which governs \eqref{egal_error}, while the approximation and generalization errors retain a similar structure for a fixed $\varepsilon>0$. Unlike the original FEONet, however, the associated constants in our setting depend implicitly on the perturbation parameter $\varepsilon$. It is noteworthy that the convergence in Theorem \ref{main_thm_whole} is not uniform with respect to $h\rightarrow0$. Indeed, this issue aligns precisely with the main theme of the reference \citep{feonet_anal}, where the authors rigorously demonstrated that both the approximation error and the generalization error depend on the condition number $\kappa(A)$ of the finite element matrix $A$, which can typically be estimated as a scale $\kappa(A)\sim h^{-2}$. This means that as $h$ becomes smaller, both the approximation and generalization errors may increase due to this adverse dependence. While this analysis was originally developed in the context of FEONet, it applies directly to eFEONet as well, since in eFEONet we solve equations with a fixed small $\varepsilon$. More precisely, the only part of the analysis in \citep{feonet_anal} where $\varepsilon$ could potentially affect the results is in the condition number estimates. If we explicitly characterize the dependency on $\varepsilon$ in these equations, then we can likewise make the $\varepsilon$
-dependence explicit in the final error estimate. In doing so, we can obtain a complete error analysis for eFEONet that incorporates both singular perturbation asymptotic analysis and the general framework from \citep{feonet_anal}, which will be addressed in the forthcoming paper.

\section{Numerical Experiments}\label{sec:experiments}
In this section, we evaluate the performance of eFEONet on three distinct types of singularly perturbed differential equations, including both ordinary and partial differential equations. For ordinary differential equations (ODEs), we examine scenarios
with and without turning points, highlighting eFEONet's adaptability to varying problems. In the case of PDEs, we assess the performance and robustness of eFEONet on two-dimensional problems defined on a square domain. In all experiments considered in this section, eFEONet employs corrector bases derived via asymptotic analysis; the derivations are provided in the appendix. Furthermore, we conduct a comparison of the experimental results with those obtained using FNO \citep{kovachki2021neural} and DeepONet \citep{lu2021learning}, the well-known operator-learning methods, and ComFNO \citep{SP03}, the neural operator model specifically designed to address the challenges of singularly perturbed differential equations. The overall experimental results are summarized in Table \ref{table:total_error}. We vary the number of training samples to assess the dependence on training data. 
\begin{table}[t]
\centering
\caption{Experimental parameters for FNO. The term ``Depth'' denotes the quantity of Fourier layers implemented within the architecture. ``LR'' designates the learning rate employed, while ``Epoch'' signifies the count of training iterations performed.}
\label{tab:fno_parameters}
\begin{tabular}{lcccc}
\toprule
\textbf{Experiment/FNO}    & \textbf{Depth}    & \textbf{LR} & \textbf{Epoch} & \textbf{Batch size} \\
\midrule
1D (no turning point)     & 4                 & 0.001       & 500            & 50                  \\ 
1D (turning point)         & 6                 & 0.001       & 500            & 50                  \\ 
2D                         & 5                 & 0.001       & 1000           & 50                  \\ 
\bottomrule
\end{tabular}
\end{table}

\begin{table}[t]
\centering
\caption{Experimental parameters for ComFNO. The term ``BlockNum'' denotes the quantity of layer blocks implemented within the architecture. ``LR'' designates the learning rate employed, while ``Epoch'' signifies the count of training iterations performed.}
\label{tab:ComFNO_parameters}
\begin{tabular}{lcccc}
\toprule
\textbf{Experiment/ComFNO} & \textbf{BlockNum} & \textbf{LR} & \textbf{Epoch} & \textbf{Batch size} \\
\midrule
1D (no turning point)      & 1                 & 0.001       & 500            & 30                  \\ 
1D (turning point)         & 2                 & 0.001       & 500            & 30                  \\ 
2D                         & 2                 & 0.001       & 1000           & 20                  \\ 
\bottomrule
\end{tabular}
\end{table}

\subsection{Experimental Setup and Implementation Details}
In this section, we describe the experimental setup used to ensure fair and reproducible comparisons across all methods. 
In order to train neural networks, we need to generate random external forcing functions. Inspired by \citep{bar2019learning}, we created a random signal $f(\xx;\om)$ as a linear combination of sine functions and cosine functions. More precisely, we set
\begin{equation}
    f(x)=m_0\sin(n_0x)+m_1\cos(n_1x),
\end{equation}for 1D cases and 
\begin{equation}
    f(x,y)=m_0\sin(n_0x+n_1y)+m_1\cos(n_2x+n_3y),
\end{equation} for 2D cases where $m_i$ and $n_j$ are sampled from the standard normal distribution.
It is worth noting that even when considering different random input functions, such as those generated by Gaussian random fields, we consistently observe similar results. This robustness indicates the reliability and stability of the eFEONet approach across various input scenarios.

The high-precision numerical solutions are denoted as $u_{\varepsilon}$, while the predictions are represented as $\widehat{u}_{\varepsilon}$. High-precision numerical solutions on the Shishkin mesh (see, e.g., \citep{SP03}) are used to compute the corresponding outputs $u_{\varepsilon}$, which serve as the ground truth during training. The training dataset comprises 900 load vectors derived from functions $f$ sampled independently. These functions are discretized using a resolution of 201 points for 1D cases and a $51 \times 51$ grid for 2D cases. Additionally, for all ODE experiments, the input-output resolution is set to 201, ensuring consistency across the comparative evaluations of FNO, ComFNO, DeepONet, and our method. In 2D PDE experiments, the resolution is fixed at 51 for $\varepsilon=10^{-3}$ and for $\varepsilon=10^{-4}$ in the rectangular domain.

We employ a neural network architecture consisting of six convolutional layers with Swish activation functions, followed by a fully connected layer. {For 1D problems, the convolutional layers are implemented with standard one-dimensional convolutions (Conv1D), whereas for 2D problems they are implemented with two-dimensional convolutions (Conv2D).} The eFEONet was trained using the L-BFGS optimizer with the following hyperparameters:
\begin{itemize}
    \item Maximal number of iterations per optimization step: $100$,
    \item Learning rate: $0.1$,
    \item History size: $100$.
\end{itemize}
All experiments were conducted on an Intel Xeon Gold 6226R CPU and an NVIDIA RTX A6000 GPU (48GB). For the 1D problems, the training dataset for FNO, ComFNO, and DeepONet includes $900\times 201$ tuples $(f, u)$, while the 2D scenarios encompass $900 \times 51 \times 51$ tuples $(f, u)$ as described in the paper \citep{SP03}. In all experiments, we employed the mean squared error (MSE) loss. For FNO and ComFNO, we used the Adam optimizer for all minimization problems, accompanied by the consistent utilization of the GELU activation function. Further details concerning the remaining parameters for our result can be found in Table \ref{tab:fno_parameters} and Table \ref{tab:ComFNO_parameters}. For DeepONet, we employed the standard branch-trunk architecture. Both the branch and trunk networks consist of three hidden layers with 100 neurons and tanh activation functions. The model was optimized using full batch and the Adam optimizer with a learning rate of $10^{-4}$ for $100,000$ training epochs.

\begin{table}[t]
  \caption{Mean relative $L^2$ test errors \textbf{($\boldsymbol{\times 10^{-3}}$)} for FNO, ComFNO, DeepONet, and eFEONet by varying the number of training input-output data pairs. Here, we set $\varepsilon=10^{-3}$ for all experiments.}
  \label{table:total_error}
  \centering
  \resizebox{\textwidth}{!}{
\begin{tabular}{ccccccccccccc}
\hline
\multirow{3}{*}{\textbf{Model}} & \multicolumn{4}{c}{Exp1. ODE w/ boundary layer} & \multicolumn{4}{c}{Exp2. ODE w/ interior layer} & \multicolumn{4}{c}{Exp3. PDE on square} \\ \cmidrule(lr){2-5} \cmidrule(lr){6-9} \cmidrule(lr){10-13}
 & \multicolumn{4}{c}{\# of training data} & \multicolumn{4}{c}{\# of training data} & \multicolumn{4}{c}{\# of training data}\\
 & 900 & 90 & 9 & None & 900 & 90 & 9 & None & 900 & 90 & 9 & None \\ \hline
FNO & 36.0 & 68.3 & 382 & - & 84.2 & 153 & 961 & - & 10.3 & 1e+03 & 1e+05 & -\\
ComFNO & 3.88 & 51.1 & 347 & - & 8.21 & 126 & 876 & - & 15.1 & 1320 & 1e+05 & -\\
DeepONet & 23.9 & 101 & 286 & - & 7.40 & 6.80 & 240 & - & 2300 & 1780 & 1590 & -\\
eFEONet (Ours) & \textbf{0.01} & \textbf{0.07} & \textbf{0.03} & \textbf{0.06} & \textbf{1.79} & \textbf{1.99} & \textbf{4.23} & \textbf{3.17} & \textbf{2.26} & \textbf{1.83} & \textbf{5.38} & \textbf{8.53} \\ \hline
\end{tabular}
}
\end{table}
\subsection{Ordinary Differential Equations with Boundary Layer}

\begin{table}[t]
    \centering
    \caption{Mean relative $L^2$ test errors\textbf{($\boldsymbol{\times 10^{-3}}$)} for FNO, ComFNO, and eFEONet across different values of $\varepsilon$ for ODEs with boundary layers.  FNO and ComFNO are trained with 900 samples, whereas eFEONet uses no pre-computed training data.}
    \label{table:case1_vary_eps}
    \resizebox{0.8\textwidth}{!}{%
        \begin{tabular}{ccccc}
            \hline
            \multirow{2}{*}{\textbf{Model}} & \multicolumn{4}{c}{Varying $\varepsilon$} \\
            & $\varepsilon = 10^{-3}$ & $\varepsilon = 10^{-4}$ & $\varepsilon = 10^{-5}$ & $\varepsilon = 10^{-6}$ \\ \hline
            FNO\;{\footnotesize(w/ 900 train data)}   & 36 & 36.8 & 36.9& 36.9 \\
            ComFNO\;{\footnotesize(w/ 900 train data)} & 3.88& 5.7 & 7.60& 5.66 \\
            Ours (eFEONet)\;{\footnotesize(w/o train data)} & \textbf{0.07} & \textbf{0.03} & \textbf{0.07} & \textbf{0.03} \\ \hline
        \end{tabular}
    }
\end{table}

\begin{figure}[t]
    \centering
    \includegraphics[width=0.7\textwidth]{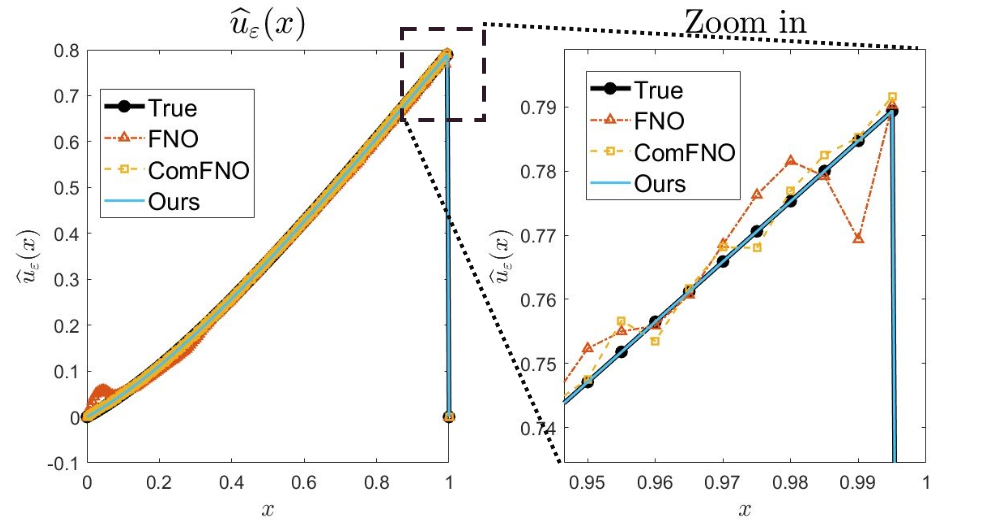}
    \caption{Comparison of predicted solutions $\widehat{u}_{\varepsilon}$
    using FNO, ComFNO, and eFEONet with $\varepsilon=10^{-4}$ for the boundary layer problem. The external force input function is given by $f(x)=1.81 \sin(1.68x)+0.09\cos(-1.78x)$.}
    \label{fig:case1_comparison}
\end{figure}

We begin with the following problem:
\begin{equation}\label{eq:case1}
    \begin{aligned}
           -\varepsilon u''_{\varepsilon} +(x+1)u'_{\varepsilon} &=f(x), \quad x\in(0,1),\\
           u_{\varepsilon}(0)&=u_{\varepsilon}(1)=0.
       \end{aligned}
\end{equation}
As shown in Figure \ref{fig:case1_comparison}, the solution exhibits an exponential boundary layer near $x=1$, making it an excellent test case for evaluating the ability of eFEONet to capture sharp boundary layers effectively. To address this challenge for \eqref{eq:case1}, eFEONet utilizes the corrector $\phi_0(x) = \exp(-2 (1-x)/\varepsilon)$ to capture the boundary layer more effectively.

As shown in the second column of Table \ref{table:total_error}, when sufficient training data is available, both FNO and ComFNO achieved reasonable accuracy, but our eFEONet outperforms them. Moreover, as the amount of training data decreases, the error for ComFNO increases significantly, whereas eFEONet maintains higher accuracy even with limited data. Table \ref{table:case1_vary_eps} presents the relative $L^2$ test errors for FNO, ComFNO, and eFEONet across different values of $\varepsilon$. The results demonstrate that eFEONet consistently outperforms the benchmark models, achieving significantly lower errors even without using any training data. Note that the error trends for FNO and ComFNO remain relatively stable across different $\varepsilon$ values, but eFEONet maintains even higher accuracy across all tested cases, demonstrating its effectiveness in capturing boundary layer phenomena without requiring extensive training datasets.

Figure \ref{fig:case1_comparison} further compares the predicted solution $\widehat{u}_{\varepsilon}$ for one of the test samples using FNO, ComFNO, and eFEONet with $\varepsilon=10^{-4}$. FNO shows substantial errors, particularly near the boundary layer, while ComFNO achieves relatively better accuracy but struggles to fully resolve the sharp transitions. In contrast, eFEONet, leveraging the corrector function as an additional basis function, achieves the highest accuracy, effectively capturing the boundary layer with minimal error.

\begin{figure}[t]
\centering
    \includegraphics[width=\linewidth]{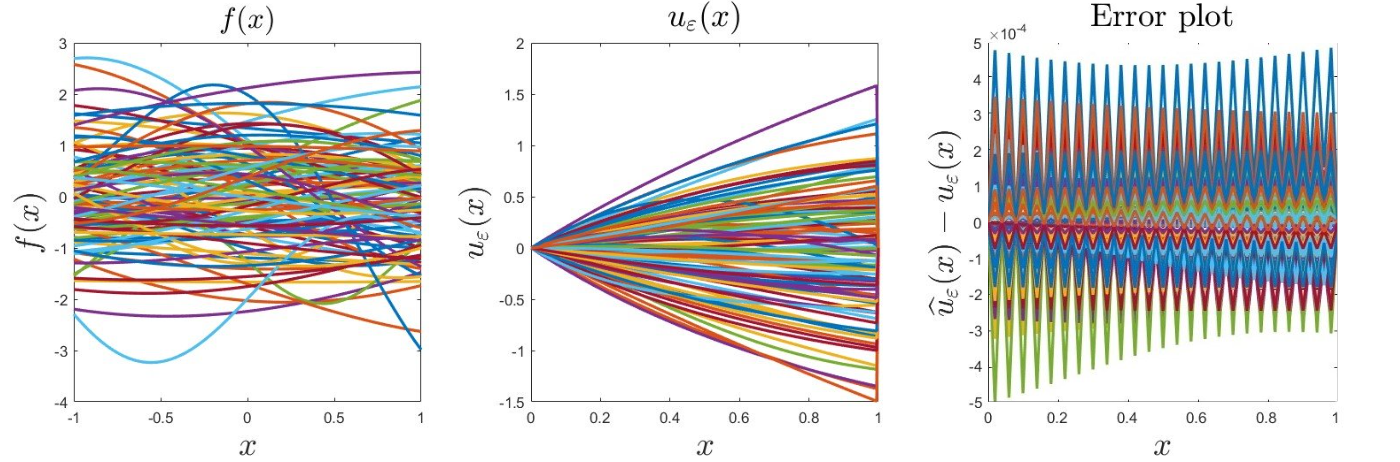} 
    \caption{Visualization of 100 input functions \( f \) (left), corresponding reference solutions (middle), and error plots (right) for boundary layer problem \eqref{eq:case1} with \( \varepsilon=10^{-5}\).}\label{appendix:case1_100}
\end{figure}

Figure \ref{appendix:case1_100} illustrates the performance of eFEONet for $\varepsilon=10^{-5}$ using 100 test samples of the forcing function $f$ at a resolution of 51. The figure on the left shows the input function \( f \), the middle figure shows the ground truth corresponding to the 100 test \( f \) samples, and the figure on the right shows the residuals produced by our method for these 100 \( f \) samples.

\subsection{Ordinary Differential Equations with Interior Layer}
We consider the following ordinary differential equation with a turning point at $x=0$:
\begin{equation}\label{eq:case2}
    \begin{aligned}
           -\varepsilon u''_{\varepsilon} -x u'_{\varepsilon} &=f(x), \quad x\in(-1,1),\\
           u_{\varepsilon}(-1)&=u_{\varepsilon}(1)=0,
   \end{aligned}    
\end{equation}   
with the corrector function $\phi_0(x)={\rm erf}(\sqrt{1/(2\varepsilon)}x)$ where ${\rm erf}(\cdot)$ is the Gauss error function. As shown in the third column of Table \ref{table:total_error}, eFEONet achieves better accuracy than both FNO and ComFNO, with a larger performance gap emerging as the number of training samples decreases. This highlights the robustness of eFEONet in data-scarce scenarios. Table \ref{table:case2_vary_eps} shows the relative $L^2$ test errors for FNO, ComFNO, and eFEONet across different values of $\varepsilon$ for ODEs with interior layers. The results demonstrate that eFEONet consistently achieves superior accuracy compared to FNO and ComFNO, even in the absence of training data. Notably, as $\varepsilon$ decreases, the performance gap between eFEONet and the benchmark models significantly widens, indicating eFEONet's ability to accurately capture sharp interior layers. 

\begin{figure}[t!]
   \centering
   \includegraphics[width=0.8\linewidth]{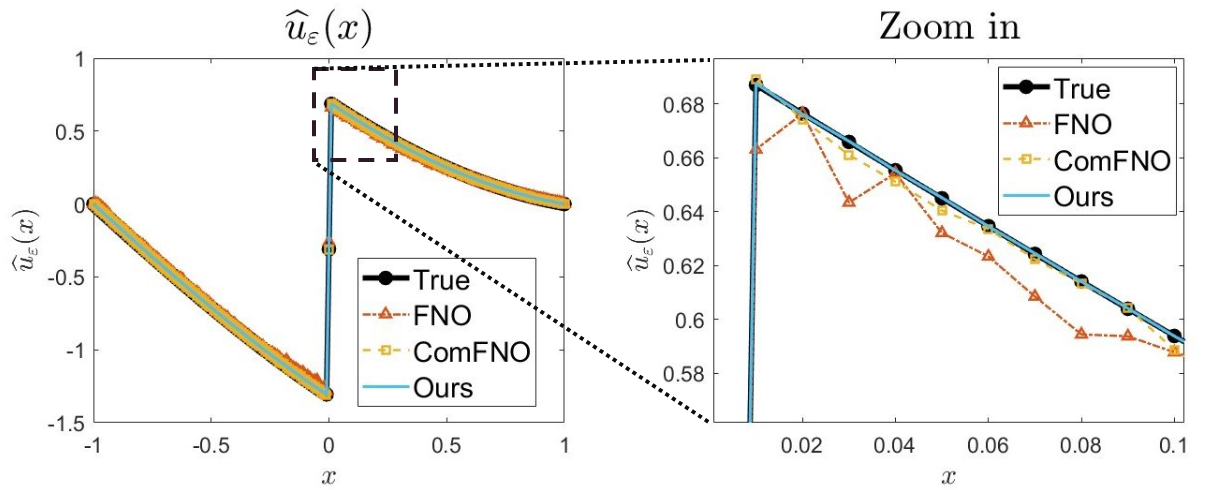}
   \caption{Comparison of predicted solutions $\widehat{u}_{\varepsilon}$
   using FNO, ComFNO, and eFEONet with $\varepsilon=10^{-8}$. The external forcing input is given by $f(x)=x(-0.58\sin(0.44x)+1.61\cos(1.05x)$).}
   \label{fig:case2_comparison}
\end{figure}

\begin{figure}[tp]
\centering
    \includegraphics[width=\linewidth]{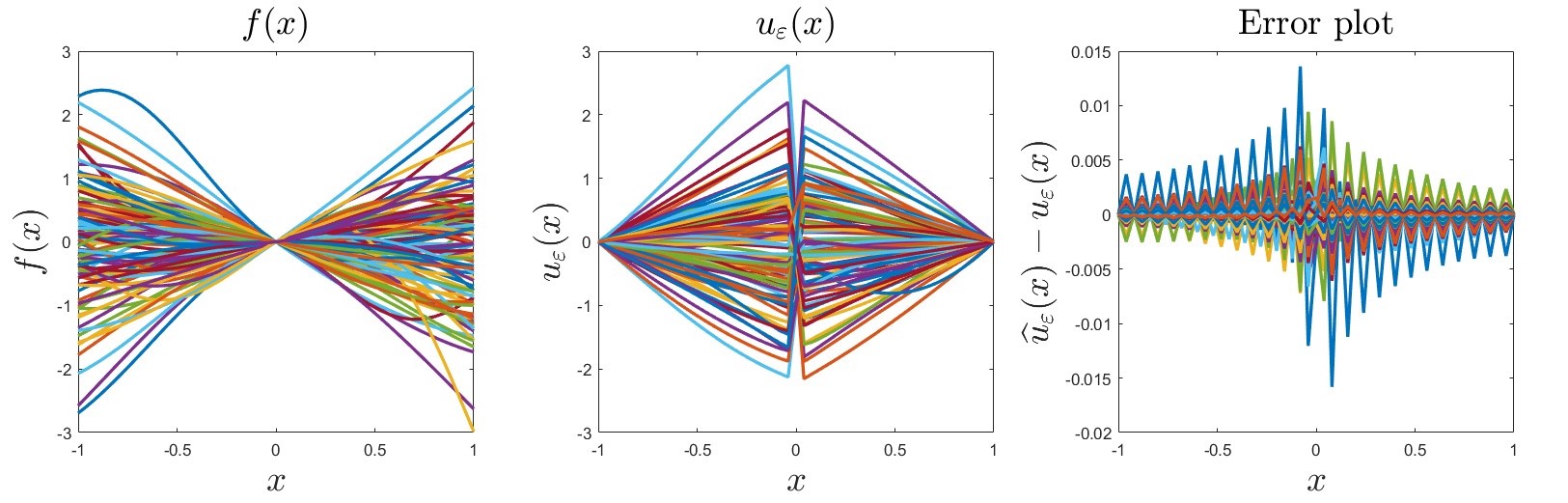} 
    \caption{Visualization of 100 input functions \( f \) (left), corresponding reference solutions (middle), and error plots (right) for interior layer problem \eqref{eq:case2} with \( \varepsilon=10^{-5}\).}\label{appendix:case2_100}
\end{figure}

Figure \ref{fig:case2_comparison} compares the predicted solutions $\widehat{u}_\varepsilon$ for two test samples using FNO, ComFNO, and eFEONet with $\varepsilon=10^{-8}$. Notably, eFEONet demonstrates superior accuracy, particularly around the singular region near the turning point at $x=0$. This result underscores the capability of eFEONet to effectively handle the challenges posed by singularities and turning points in differential equations, delivering reliable predictions even in complex scenarios.  Figure~\ref{appendix:case2_100} shows the input function \( f \) (left), the ground truth corresponding to the 100 test \( f \) samples (middle), and the residuals produced by our method for these 100 \( f \) samples (right) at a resolution of 51.

\begin{table}[t]
\caption{Mean relative $L^2$ test errors\textbf{($\boldsymbol{\times 10^{-3}}$)} for FNO, ComFNO, and eFEONet across different values of $\varepsilon$ for ODEs with interior layers. The results highlight the performance of each model when trained with 900 data samples (FNO, ComFNO) and without training data (eFEONet).}
\label{table:case2_vary_eps}
\centering
\resizebox{0.86\textwidth}{!}{
\begin{tabular}{ccccc}
\hline
\multirow{2}{*}{\textbf{Model}} & \multicolumn{4}{c}{Varying $\varepsilon$} \\
& $\varepsilon=10^{-3}$ & $\varepsilon=10^{-4}$ & $\varepsilon=10^{-5}$ & $\varepsilon=10^{-6}$ \\ \hline
FNO{\footnotesize w/ 900 train data} & 84.2 & 86.9 & 81.5 & 86.9 \\
ComFNO{\footnotesize w/ 900 train data} & 8.21 & 8.97 & 19.6 & 15.5 \\
Ours(eFEONet){\footnotesize w/o train data} & \textbf{3.17} & \textbf{5.21} & \textbf{0.66} & \textbf{0.19} \\ \hline
\end{tabular}
}
\end{table}
\subsection{Partial Differential Equations on Square}
\begin{figure}[t]
    \centering
    \includegraphics[width=\linewidth]{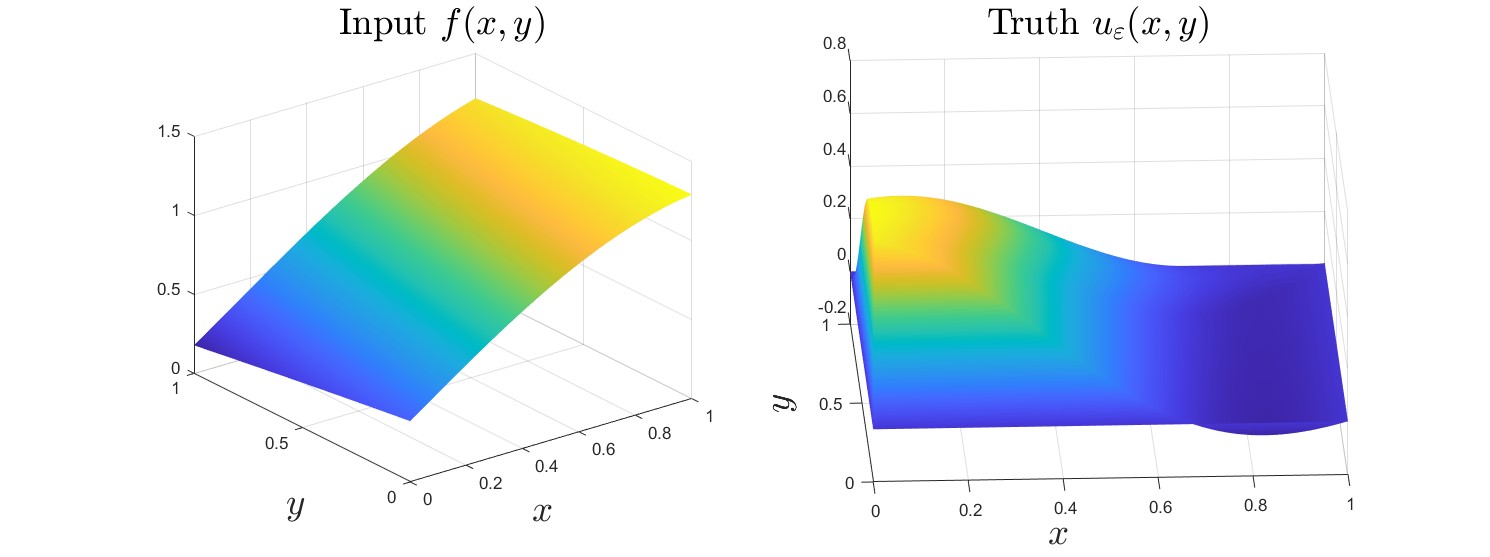}
    \caption{Solution profiles for the PDE problem on a square domain.}
    \label{fig:case3_solution profile}
\end{figure}

\begin{figure}[t]
   \centering
   \includegraphics[width=0.9\linewidth]{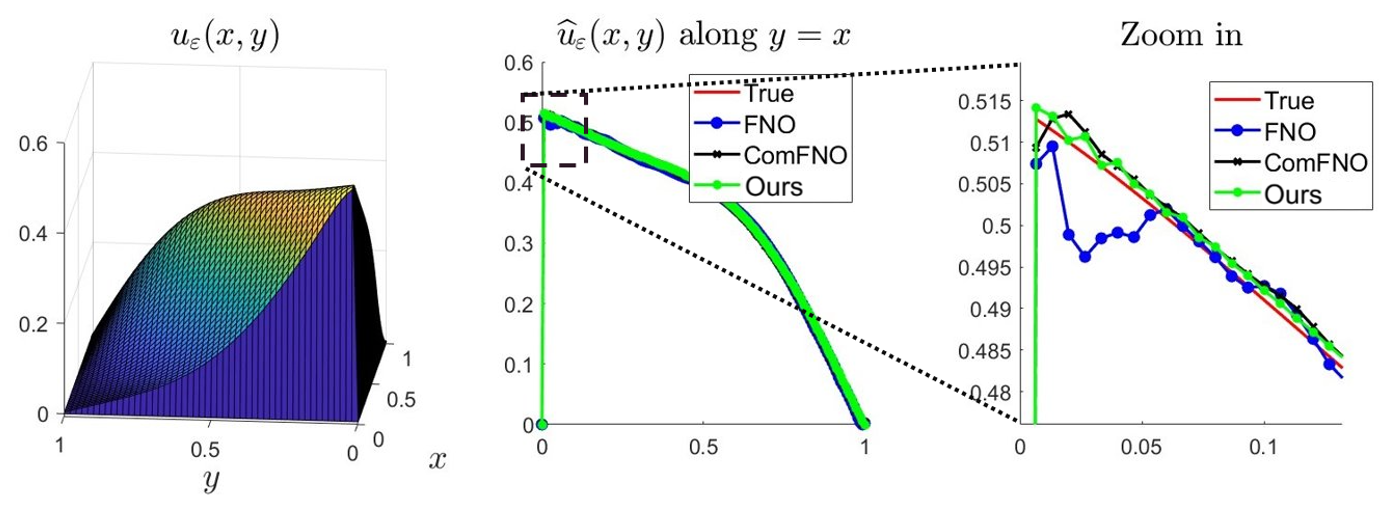}
   \caption{Comparison of the reference solution $u_{\varepsilon}(x,y)$ (left) and the predicted solutions $\widehat{u}_{\varepsilon}$ along the diagonal $y = x$ plane for ComFNO and eFEONet (middle and right) with $\varepsilon=10^{-4}$. The results highlight the superior accuracy of eFEONet in capturing sharp boundary layers along $x=0$, whereas ComFNO exhibits noticeable errors near the boundary regions.}
   \label{fig:case3_error comaparison}
\end{figure}

\begin{figure}[htp]
\centering
    \includegraphics[width=\linewidth]{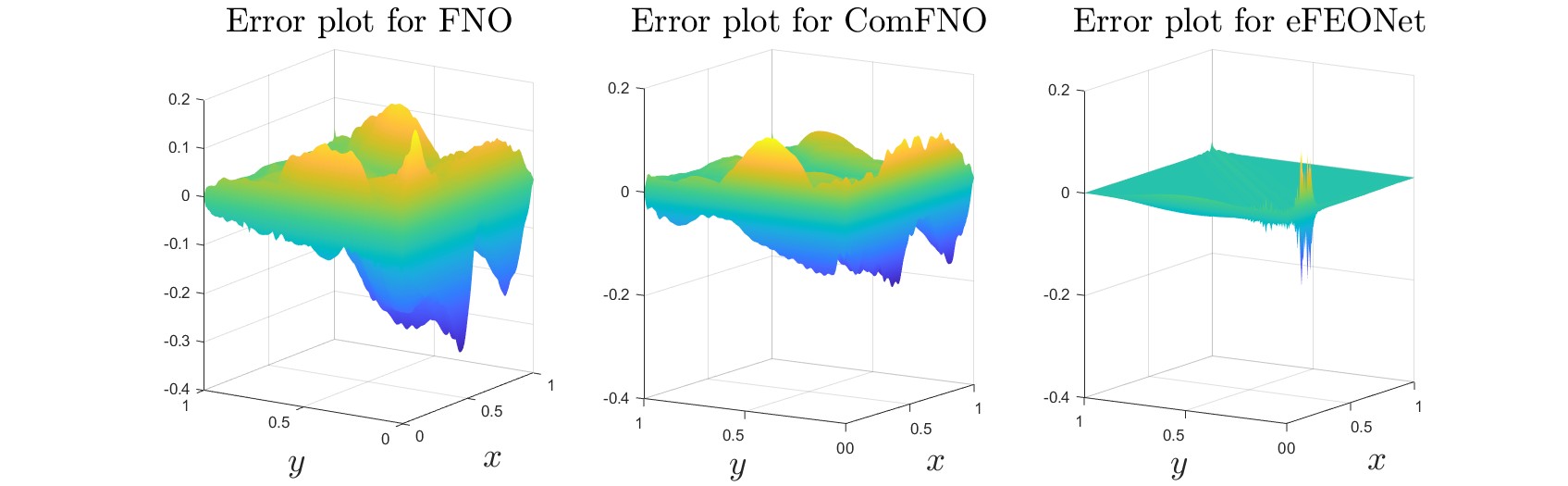} 
    \caption{Error visualizations for 100 test forcing functions using FNO (left), ComFNO (middle), and eFEONet (right) for the PDE \eqref{eq:case3} on square with $\varepsilon = 10^{-4}$ and input–output resolution $51$.}\label{appendix:case3_100}
\end{figure}

For a boundary-value problem of an elliptic PDE on the spatial domain $D=[0,1]^2$, we consider
\begin{equation}\label{eq:case3}
    \begin{aligned}
        -\varepsilon \Delta u_{\varepsilon}-(1,1)\cdot \nabla u_{\varepsilon}&=f(x,y) \,\,\,\,\, \text{in}~ D,\\
        u_{\varepsilon}(x,y)&=0\quad\,\qquad\text{on}~ \partial D,
    \end{aligned}
\end{equation}
where the solution exhibits a boundary layers along the edges at $x=0$ and $y=0$, as illustrated in Figure \ref{fig:case3_solution profile}. For this PDE problem, the asymptotic expansion of \( u(x,y) \) can be expressed as:
\[
    u(x,y) = u_0(x,y) - u_0(0,y)e^{-x/\varepsilon} - u_0(x,0)e^{-y/\varepsilon}+ u_0(0,0) e^{(-x-y) /\varepsilon}.
\]
As shown in the fourth column of Table \ref{table:total_error}, the accuracy gap between eFEONet and benchmark models becomes even more pronounced for this problem. This highlights the capability of eFEONet to effectively resolve boundary layers in complex spatial domains.
Furthermore, Figure~\ref{fig:case3_error comaparison} illustrates that ComFNO exhibits relatively large errors near the boundary layers, whereas eFEONet maintains uniformly low error levels across the domain. This highlights the robustness of eFEONet in accurately resolving sharp boundary-layer structures. We overlaid the error plots for 100 random inputs of $f$ in Figure~\ref{appendix:case3_100}, calculated with $\varepsilon = 10^{-4}$.

\subsection{Comparison with the original FEONet}
\begin{figure}[tp]
    \centering
    \includegraphics[width=0.85\linewidth]{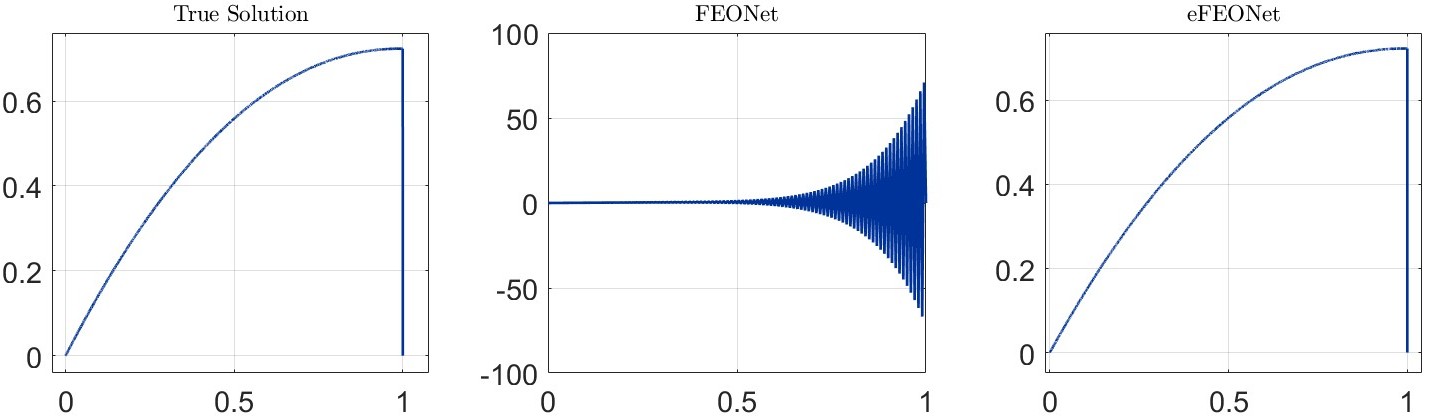}
    \caption{Visualization of the true solution and solution plots from FEONet and eFEONet of the boundary layer problem when $\varepsilon=10^{-5}$.}\label{fig: section3.5}
\end{figure}

The original FEONet is limited by its reliance on piecewise polynomial basis functions, which are insufficient for resolving stiff behaviors. As a result, it offers no significant performance advantage over the standard FEM (see the middle plot in Figure \ref{fig: section3.5}). To address this, we proposed eFEONet, which enriches the basis functions with exponential and error functions derived from the asymptotic analysis (see Appendix~\ref{sec:app_deriv} for details), ensuring they align with the mathematical structure of the problem's stiff behavior. This enhancement enables eFEONet to outperform the original FEONet, yielding substantially reduced errors across both boundary and interior layers, as shown in Table \ref{tab:feonet_comparison}.
\begin{table}
\centering
\caption{Comparison of FEONet and eFEONet for $\varepsilon = 10^{-5}$. 
Errors are reported for the boundary layer and interior layer regions.}
\begin{tabular}{lcc}
\hline
\textbf{Model} & \textbf{Boundary Layer} & \textbf{Interior Layer} \\ \hline
FEONet & 3.04 & 0.0222 \\
Ours (eFEONet) & \textbf{7.0e-05} & \textbf{6.6e-04} \\
\hline
\end{tabular}
\label{tab:feonet_comparison}
\end{table}
\section{Conclusion}

In this paper, we introduced eFEONet, designed for singularly perturbed differential equations. By integrating boundary layer theory into the finite element framework, eFEONet captures sharp transitions using theory-guided basis functions, eliminating the need for extensive training datasets.
Experimental results demonstrate the robustness of eFEONet across various PDEs with boundary, interior layers, and corner layer problems. Compared to FNO and ComFNO, eFEONet consistently achieves superior accuracy, particularly in data-scarce scenarios. Additionally, our method is supported by convergence analysis, validating its reliability.
Despite its strong performance, certain limitations remain. First, the choice of parameters, such as the number of basis functions and network hyperparameters, significantly affects the learning dynamics and overall performance of eFEONet. A systematic analysis of these parameters is still an open research question. Second, while our study presents a unique method for solving singularly perturbed problems with boundary and interior layers using minimal or even no training data, future research should extend eFEONet to handle more challenging problems, such as corner singularities and other intricate geometrical effects.

\section*{Acknowledgments}
Y. Hong was supported by 
Basic Science Research Program through the National Research Foundation of Korea (NRF) funded by the Ministry of Education (NRF-2021R1A2C1093579) and by the Korea government(MSIT) (RS-2023-00219980).
S. Ko was supported by National Research Foundation of Korea Grant funded by the Korean Government (RS-2023-00212227).
Jae Yong Lee was supported by Institute for Information \& Communications Technology Planning \& Evaluation (IITP) through the Korea government
(MSIT) under Grant No. 2021-0-01341 (Artificial Intelligence Graduate School Program (Chung-Ang University)), and by the National Research Foundation of Korea(NRF) grant funded by the Korea government(MSIT) (RS-2025-02303239).	

\appendix
\appendixpage
\addappheadtotoc
\section{Derivation of Corrector Basis Functions}\label{sec:app_deriv}
We have focused on the numerical treatment of the following singularly perturbed convection-dominated problem
\begin{align*}
    -\varepsilon \Delta u_\varepsilon - \boldsymbol{b} \cdot \nabla u_\varepsilon + c u_\varepsilon &= f \quad \text{in } D,\\
    u_\varepsilon &= 0 \quad \text{on } \partial D,
\end{align*}
where $0 < \varepsilon \ll 1$, and $ \boldsymbol{b} = \boldsymbol{b}(\mathbf{x})$, $c = c(\mathbf{x})$ and $f = f(\mathbf{x})$ are given smooth functions defined over the domain \( D \). This formulation represents a general convection-diffusion-reaction equation with singular perturbation. For this problem, we considered both 1D and 2D settings, addressing critical challenges such as boundary layers and interior layers that arise due to the small parameter \( \varepsilon>0 \). 
From this point onward, our analysis follows the singular perturbation analysis stated in \citep{BLA_3}. The theoretical foundations and techniques presented here are based on this approach, providing a rigorous framework for handling boundary and interior layers in singularly perturbed problems. For further details on related studies and extensions, we refer the reader to \citep{BLA_3}.

{\bf Boundary layer case.} We simplify the analysis and explanation by considering a one-dimensional paradigm problem. The 1D problem is defined as
\begin{align*}
    -\varepsilon u_\varepsilon'' - u_\varepsilon' &= f \quad {\rm{in}}\,\,(0,1),\\
    u_\varepsilon(0) &= u_\varepsilon(1) = 0.
\end{align*} 
This 1D model provides a clear framework for understanding boundary layer phenomena and allows us to systematically develop the necessary mathematical and computational tools before extending the approach to higher dimensions.
The corresponding limit problem is obtained by formally setting \( \varepsilon = 0 \):
\begin{align*}
    - u_0' &= f \quad {\rm{in}}\,\,(0,1),\\
    u_0(1) &= 0.
\end{align*}
Treating this as a transport equation, we supplement the limit problem with the inflow boundary condition at \( x = 1 \), namely
\[
u_0(1) = 0.
\]
Solving this equation with the given condition yields
\[
u_0 = - \int_x^1 f(s) \, {\rm{d}}s.
\]
At this stage, the choice of the inflow boundary condition \( u_0(1) = 0 \) is an assumption motivated by the structure of the transport equation. 
To address the boundary layer near \( x = 0 \), we introduce a stretched variable \( \bar{x} = x / \varepsilon^\alpha \), with \( \alpha > 0 \). Substituting \( \bar{x} \) into the original problem with \( f = 0 \), we derive
\[
-\varepsilon^{1-2\alpha} \frac{{\rm{d}}^2 u_\varepsilon}{{\rm{d}}\bar{x}^2} - \varepsilon^{-\alpha} \frac{{\rm{d}} u_\varepsilon}{{\rm{d}}\bar{x}} = 0.
\]
Here, \( f \) is omitted because it is accounted for in the inviscid equation \( -u_0' = f \). To define a corrector from this equation, we observe that the corrector must balance the difference between \( u_\varepsilon \) and \( u_0 \) at \( x = 0 \) and decay rapidly as \( x \) moves away from 0.
By setting \( 1 - 2\alpha = -\alpha \), we find \( \alpha = 1 \), resulting in the following boundary layer equation
\[
-\frac{{\rm{d}}^2 \bar{\theta}_\varepsilon}{{\rm{d}}\bar{x}^2} - \frac{{\rm{d}} \bar{\theta}_\varepsilon}{{\rm{d}}\bar{x}} = 0.
\]
The boundary conditions for this equation are
\[
\bar{\theta}_\varepsilon(0) = -u_0(0), \quad \bar{\theta}_\varepsilon \to 0 \quad \text{as } \bar{x} \to \infty.
\]
The explicit solution for \( \bar{\theta}_\varepsilon \), the approximate corrector, is given as
\[
\bar{\theta}_\varepsilon = -u_0(0) e^{-\bar{x}} = -u_0(0) e^{-x / \varepsilon}.
\]
As discussed earlier, we want to add this boundary layer function into our finite element ansatz space. However, note that this boundary layer function does not satisfy the appropriate boundary conditions. This is easily handled by introducing the boundary layer basis function of the form
\[
    \phi_0(x)=e^{-x/\varepsilon}+(1-e^{-1/\varepsilon})x+1.
\]

{\bf Interior layer case.} For convection-diffusion equations with an interior layer, we consider the problem
\begin{align*}
    -\varepsilon \ u_\varepsilon'' - b(x) u_\varepsilon' &= f \quad  {\rm{in}}\,\, (-1, 1),\\
    u_\varepsilon(-1) &= u_\varepsilon(1) = 0,
\end{align*}
where \( b(x) \) satisfies \( b < 0 \) for \( x < 0 \), \( b(0) = 0 \), \( b > 0 \) for \( x > 0 \), and \( b'(x) > 0 \). The turning point at \( x = 0 \) introduces an interior layer due to the change in sign of \( b(x) \), where characteristics collide.
For the formal limit problem, setting \( \varepsilon = 0 \) leads to
\[
-b(x) u_0' = f.
\]
However, it may not be well-defined at \( x = 0 \) since \( b(0) = 0 \). Therefore, we split the solution into left and right parts, \( u_0^l \) and \( u_0^r \), corresponding to \( x < 0 \) and \( x > 0 \), respectively
\[
-b(x) (u_0^l)' = f \quad \text{for } x < 0,\quad{\rm{and}} \quad -b(x) (u_0^r)' = f \quad \text{for } x > 0.
\]
The inflow boundary conditions are then supplemented as
\[
u_0^l(-1) = 0, \quad u_0^r(1) = 0.
\]
The discrepancy at \( x = 0 \) between \( u_0^l \) and \( u_0^r \) produces an interior layer. If \( f(0) = 0 \), the correctors introduced below can effectively capture the sharpness of this layer. However, if \( f(0) \neq 0 \), the limit problem
\[
-b(x) u_0' = f
\]
has an inconsistency at \( x = 0 \) because \( b(0) = 0 \). This implies that \( u_0' \) diverges near \( x = 0 \), and the interior layer cannot be fully captured by standard corrector functions. To address this issue, the data may need to be adjusted to ensure compatibility, as described in related perturbation analyses.
To analyze the interior layer, we introduce the stretched variable \( \bar{x} = x / \sqrt{\varepsilon} \) and approximate \( b(x) \) as \( b(x) = b'(0)x + \frac{1}{2}b''(\xi)x^2 \approx b'(0)\sqrt{\varepsilon}\bar{x} \). Substituting these into the original equation with \( f = 0 \), we obtain the leading-order differential equation
\[
-\frac{{\rm{d}}^2\theta}{{\rm{d}}\bar{x}^2} - b'(0) \bar{x} \frac{{\rm{d}}\theta}{{\rm{d}}\bar{x}} = 0,
\]
subject to the boundary conditions
\[
\theta \to \text{constant as } \bar{x} \to \pm \infty.
\]
The solution of this equation, written explicitly, is
\[
\theta = \frac{2}{\sqrt{\pi}} \int_0^{\bar{x} \sqrt{b'(0) / 2}} e^{-\tau^2}\, {\rm{d}}\tau = \text{erf}\left( \bar{x} \sqrt{\frac{b'(0)}{2}} \right)= \text{erf}\left( x \sqrt{\frac{b'(0)}{2\varepsilon}} \right),
\]
where \( \text{erf} \) denotes a Gaussian error function. This serves as a corrector for the interior layer.

\bibliography{references}
\bibliographystyle{abbrv}
\end{document}